\documentclass[12pt]{article}
\usepackage{amsmath,amssymb,amstext,dsfont,fancyvrb,float,fontenc,graphicx,subfigure,theorem,hyperref}

\usepackage{tikz,everypage}
\usepackage[letterpaper]{geometry}
\newfont{\footsc}{cmcsc10 at 8truept}
\newfont{\footbf}{cmbx10 at 8truept}
\newfont{\footrm}{cmr10 at 10truept}

\usepackage{fancyhdr}
\usepackage{relsize}
\usepackage{sectsty}
\allsectionsfont{\larger[-1]} 

\renewcommand\paragraph{\@startsection{paragraph}{4}{\z@}
                                    {2ex \@plus.5ex \@minus.2ex}
                                    {-1em}
                                    {\normalfont\normalsize\bfseries}}

\renewcommand\subparagraph{\@startsection{subparagraph}{5}{\parindent}
                                       {2ex \@plus.5ex \@minus .2ex}
                                       {-1em}
                                      {\normalfont\normalsize\bfseries}}

\newlength{\BiblioSpacing}
\renewenvironment{thebibliography}[1]{
\begin{oldthebibliography}{#1}
\setlength{\parskip}{\BiblioSpacing}
\setlength{\itemsep}{\BiblioSpacing}
}
{
\end{oldthebibliography}
}

\usepackage[strict]{changepage}
\def\abstractname{Abstract -}   % <-----------------
\def\abstract{\begin{adjustwidth}{1cm}{1cm} \par    \footnotesize \noindent {\bf \abstractname} 
\def\endabstract{ \end{adjustwidth} \smallskip }}
{\theorembodyfont{\itshape}\newtheorem{theorem}{Theorem}[section]}
{\theorembodyfont{\itshape}}
{\theorembodyfont{\itshape}}
{\theorembodyfont{\itshape}}
{\theorembodyfont{\itshape}}
{\theorembodyfont{\rm}}
{\theorembodyfont{\rm}}
{\theorembodyfont{\rm}\newtheorem{conjecture}[theorem]{Conjecture}}
{\theorembodyfont{\rm }}

\usepackage{csquotes} 
\newcommand{\Lim}[1]{\raisebox{0.5ex}{\scalebox{0.8}{$\displaystyle \lim_{#1}\;$}}}

\title{\Large\bf Semicomplete Arithmetic Sequences, Division of Hypercubes, and the Pell Constant}
\author{\sc Z. Hoelscher}

\numberwithin{equation}{subsection}
\numberwithin{theorem}{subsection}

\begin{document}
\setcounter{page}{1}
%\date{}
%%%%%%%%%%%%%%%%%%%%%%%%%%%%%%%%%%%%%%%%%
\maketitle
\thispagestyle{fancy}

\vskip 1.5em

\begin{abstract}
In this paper we produce a few continuations of our previous work on partitions into fractions. Specifically, we study strictly increasing integer sequences $\{n_j\}$ such that there are partitions for all integers less than the floor of $G$, where $G=\frac{n_1}{j} + \frac{n_2}{j} + \cdots + \frac{n_{j-1}}{j}+\frac{n_j}{j}$, and all summands are distinct terms drawn from $\frac{n_1}{j} + \frac{n_2}{j} + \cdots + \frac{n_{j-1}}{j}+\frac{n_j}{j}$. We call such sequences \enquote{semicomplete}. We find that there are only three semicomplete arithmetic sequences. We also study sequences that give the maximum number of pieces that an $M$ dimensional hypercube can be cut into using $N-1$ hyperplanes. We find that these are semicomplete in one, two, three, and four dimensions. As an aside, we use one of our generating functions to produce what appears to be a new identity for the Pell constant, a number which is closely connected to the density of solutions to the negative Pell equation. 
\end{abstract}
 

\begin{MSC}
05A17, 11B65
\end{MSC}
\section{Introduction}

Integer partitions have a long history of study, with mathematicians from Euler to Ramanujan producing interesting results. Past mathematicians have also studied complete sequences. These are defined in the work of Brown \cite{one} as sequences $\{ f_i \}^{\infty}_{i=1}$ such that all terms are positive and every $n\in\mathbb{N}$  can be written as shown below. 
\begin{equation}
    n=\sum^{\infty}_{i=1} \alpha_i f_i \hspace{0.2 in} \alpha_i \in \{ 0, 1\}
\end{equation}

In this paper we build on our previous work, where we studied partitions of integers into fractions with constant denominators and distinct numerators drawn from sets of even or odd integers \cite{two}. We noted that for any positive integer $j$, one can write $j^2 = \sum_{n=1}^j (2n-1)$. After dividing both sides of the equation by $j$, one can represent $j$ as a sum of fractions. 
\begin{equation}
    j= \frac{1}{j}+\frac{3}{j}+\frac{5}{j}+\cdots+\frac{2j-1}{j}
\end{equation}
We then proved that when $j>2$, any $k \in \mathbb{N}, k<j$ can be written as the sum of some combination of distinct terms from the series that sums to produce $j$. We noted that this is not always possible when the numerators are even integers, prompting us to consider what other sequences allow for similar behavior. Such a generalization of the problem is considered in this paper. Suppose we have a series of the form given below:
\begin{equation}
G=\frac{n_1}{j} + \frac{n_2}{j} + \cdots + \frac{n_{j-1}}{j}+\frac{n_j}{j}
\end{equation}

We look for solutions for natural numbers $k<\lfloor G \rfloor$, where $k = \frac{a_1}{j} + \frac{a_2}{j}+\cdots +\frac{a_{m-1}}{j}+\frac{a_m}{j}$, and $a_1, \ldots, a_m$ are distinct coefficients drawn from $\{n_1, n_2,\ldots,n_j\}$. We also require that $k,j \in \mathbb{N}$, $j>2$. We say $\{ n_j \}$ is semicomplete if it enables such solutions for all $k$, where the terms of the sequence strictly increase, and all terms are positive. We refer to such sequences as semicomplete because a semicomplete sequence is not necessarily complete. For example, we prove that the cake numbers are both complete and semicomplete, yet the odd integers are only semicomplete. This is because the odd integers enable sums for all $kj$ under the restrictions above, though one cannot find integers such as 2 as the sum of distinct odd summands. This makes semicomplete sequences perhaps a bit less obvious, as one does not have to be able to find partitions for every natural number. Such questions lead naturally to the following theorem.

\begin{theorem}
\label{thm1}
The only semicomplete arithmetic sequences are $\{1, 3, 5, 7, 9\ldots \}$, \newline  $\{2, 3, 4, 5, 6 \ldots \}$, and $\{1, 2, 3, 4, 5 \ldots \}$.
\end{theorem}
We prove this theorem in Section \ref{sec2}. 

In our previous paper, we conjectured that the cake numbers and lazy caterer's sequence are both semicomplete \cite{two}. We prove this now, and generalize these sequences to higher dimensions. Let $C_M^{N}$ denote the maximum number of pieces that an $M$-dimensional hypercube can be cut into using $N-1$ hyperplanes. Let $\{C_M^{N}\}$ be the sequence of such terms for a given value of $M$ where $N$ varies. (Observe that when $M=2$ we have the lazy caterer's sequence, and when $M=3$, we have the cake numbers.) Note that we use $N-1$ for the number of cuts so that $N=1$ gives the first term of the sequence, $N=5$ gives the fifth term, et cetera, rather than $N=0$ for the first term and $N=4$ for the fifth. 
\begin{theorem}
\label{thm2}
For values of $M$ equal to one, two, three, or four, $\{C_M^N\}$ is semicomplete. \end{theorem}
We provide a proof for this in Section \ref{sec3}. We also give an interesting conjecture on these sequences.
\begin{conjecture}
$\{C_M^N\}$ is semicomplete for any value of $M\geq 1$, $M \in \mathbb{N}$.
\end{conjecture}

If true, this would imply that an infinite number of semicomplete sequences exist, hence one cannot list them all. To solve the full semicompleteness problem, one would then have to find the general necessary and sufficient conditions for any type of sequence to be semicomplete. 

Perhaps a more interesting problem would be to search for sequences that, like the odd integers, are semicomplete but not complete. We know of one example, $\{2,3,4,5,6.\dots\}$, though it would be interesting to find more. 

We note that if one does not require $k$ to be an integer, one can write the generating function for the case of odd numerators as shown below. 
\begin{equation}
    \prod_{n=1}^{j} (1+q^{2n-1}) = \sum_{kj=0}^{j^2} f_{O_j}(k,j) q^{kj}
\end{equation}
As a somewhat interesting aside, we give what appears to be a new $q$-series identity for the Pell constant, $\mathcal{P}_{Pell}$. This constant is closely connected to the product $\prod_{n=1}^{j} (1+q^{2n-1})$ as well as the density of solutions to the negative Pell equation $x^2 - Dy^2 = -1$ \cite{three}. It is currently unknown whether the Pell constant is transcendental, though it has been proven to be irrational \cite{three}.
\begin{theorem}
\label{thm3}
\begin{equation}
    F(q)=\sum_{m=0}^{j+1} \sum_{n=0}^{j+1} \frac{(-1)^{\frac{3n}{2}}q^{\frac{1}{2}(n^2-2n+2)} (-1)^{\frac{m}{2}}q^{\frac{1}{2}(m^2-2m)}}{q+1} {j+1 \brack n}_q {j+1 \brack m}_q 
\end{equation}
$\mathcal{P}_{Pell}=1- \Lim{j \rightarrow \infty} F(-\frac{1}{2}) \approx 0.58057$
\end{theorem}
We provide a proof for this in Section \ref{sec4}. We note that ${j+1 \brack n}_q$ is a Gaussian binomial coefficient, which is the $q$-analog of a binomial coefficient. Gaussian binomial coefficients are defined as shown below \cite{four}, where $m\geq n$.
\begin{equation}
{m \brack n}_q = \frac{(1-q)(1-q^{2}) \cdots (1-q^{m})}{(1-q)(1-q^2) \cdots (1-q^n)(1-q)(1-q^2) \cdots (1-q^{m-n})}
\end{equation}
One should note that when $n>m$, ${m \brack n}_q$ is defined as zero. We can then restate our conjecture from our previous paper \cite{two} in terms of $F(q)$, where we note that $F(q)$ is a polynomial in $q$ when fully simplified. 
\begin{conjecture}
Let $G(q)$ contain only the terms from $F(q)$ where the exponent of $q$ is $kj$, $k,j \in \mathbb{N}, k<j, j>2$. The coefficients of $G(q)$ are always either unimodal or bimodal. 
\end{conjecture}

\section{Proof for Theorem \ref{thm1}}
\label{sec2}
As we require that we must have solutions for all $k<\lfloor G \rfloor$ where $j>2$, $k,j\in \mathbb{N}$, we can prove that a sequence does not work by finding a value of $k$ that has no solution for some value of $j$. While somewhat tedious, this process works well to eliminate all sequences other than the three semicomplete arithmetic sequences. We can then readily show these three are semicomplete. 

\begin{proof}
When $j=3$, we have a series of the form
\begin{equation}
    G=\frac{a}{3}+\frac{a+b}{3}+\frac{a+2b}{3}
\end{equation}
where $a, b \in \mathbb{N}$. Note that $a$ is the first term in the sequence of numerators and $b$ is the difference between consecutive terms in that sequence. 

\noindent Let $a=1$:
\begin{equation}
    G=\frac{1}{3}+\frac{1+b}{3}+\frac{1+2b}{3}
\end{equation}

If $b=0$, we have the sequence $\{1, 1, 1, 1, \ldots \}$ as numerators. In this case $G$ is always exactly one, hence there are no natural numbers $k<\lfloor G \rfloor$. We thus see that $\{1, 1, 1, 1, \ldots \}$ is not semicomplete. If $b=1$, we have the consecutive integers as numerators. It is easy to see that any integer less than the sum of the first $n$ consecutive integers can be found as the sum of some combination of terms from the first $n$ consecutive integers, hence this sequence is semicomplete. 

If $b=2$ our numerators are the odd integers, which we have previously proven to be semicomplete \cite{two}. If $b>2$ all terms in the series but $\frac{1}{3}$ are greater than one, hence it is impossible to find a combination of terms that sum to $k=1$. This then proves that any sequence with $a=1$ and $b>2$ is not semicomplete. 

We now consider sequences where $a>1$. If $b=0$ and $a=2$ there is no solution for $k=1$, as $1$ is not an integral multiple of $\frac{2}{3}$. If $a=3$ and $b=0$ there are solutions for all $k$ when $j=3$, but not when $j=4$, hence $\{3, 3, 3, 3, \ldots \}$ is not semicomplete. If $a>3$ and $b=0$ all terms in the series 
\begin{equation}
    G=\frac{a}{3}+\frac{a+b}{3}+\frac{a+2b}{3}
\end{equation}
are greater than one, hence there is no combination for $k=1$. All sequences $\{a, a, a, a, \ldots \}$, where $a>3$ are thus incomplete.  

If $a=2$ and $b=1$ we have the sequence $\{2, 3, 4, 5, \ldots \}$. One can see that any integer except $1$ and $(2+3+4+5+\cdots+(j+1))-1$ can be found as the sum of some combination of terms from $\{2, 3, 4, 5, \ldots \}$. We know that $\frac{1}{j}$ and $\frac{(2+3+4+5+\cdots+(j+1))-1}{j}$ cannot be integers when $j>2$, hence every integer $k<\lfloor G \rfloor$ can be found as the sum of some combination of terms. We know $\frac{(2+3+4+5+\cdots+(j+1))-1}{j}$ cannot be an integer through the argument shown below. 

\begin{equation}
    G=\frac{(j+1)(j+2)}{2j}-\frac{1}{j}=\frac{j+3}{2}
\end{equation}
\begin{equation}
    \frac{(2+3+4+5+\cdots+(j+1))-1}{j}=\frac{j+3}{2}-\frac{1}{j}=\frac{j}{2}-\frac{1}{j}+\frac{3}{2}
\end{equation}
We know that $\frac{j}{2}-\frac{1}{j}+\frac{3}{2}$ cannot be an integer when $j>2$, hence $\frac{(2+3+4+5+\cdots+(j+1))-1}{j}$ cannot be an integer when $j>2$.

If $a=2$ and $b>1$ then there is no combination for $k=1$ when $j=3$, as every term in the series but $\frac{2}{3}$ is greater than $1$. Such sequences are thus incomplete. If $a=3$ and $b=1$ or $b=2$, there is no combination for $k=2$ when $j=3$. If $a=3$ and $b=3$ there is no combination for $k=1$ when $j=4$. If $a=3$ and $b>3$ there is no combination for $k=2$ when $j=3$. If $a>3$ then all terms in the series are greater than $1$ when $j=3$, regardless of the value of $b$, hence there is no combination for $k=1$. Such sequences are thus incomplete. 
\end{proof}

\section{Proof for Theorem \ref{thm2}}
\label{sec3}
\begin{proof}
\textbf{In one dimension:}
It is easy to see that in this case we are cutting a line into pieces with points. This then results in $\{C_1^{N}\}=\{1,2,3,\ldots\}$, which is the sequence of consecutive integers. We know this sequence is semicomplete. 
 
\noindent \textbf{In two dimensions:}
\begin{equation}
    \{C_2^{N}\}=\bigg\{\frac{N^2 - N + 2}{2}\bigg\}
\end{equation}
The sum of the first $t$ terms can be found through induction. 
\begin{equation}
    \sum_{N=1}^{t}\frac{N^2 - N + 2}{2} = \frac{1}{6}t(t^2 + 5)
\end{equation}
We see that when $t$ is sufficiently large, the $(t+1)^{TH}$ term is less the sum of the first $t$ terms. We have 
\begin{equation}
    \frac{(t+1)^2-(t+1)+2}{2} < \frac{1}{6}t(t^2 + 5) 
\end{equation}
\noindent where $3 < t$.

Suppose one takes a set of the first $t$ terms and one sees that any integer less than the sum of those $t$ terms can be found as the sum of some combination of terms from that set. Now take an integer $I$ that satisfies the following inequality:
\begin{equation}
   \sum_{N=1}^{t}\frac{N^2 - N + 2}{2} < I < \sum_{N=1}^{t+1}\frac{N^2 - N + 2}{2}
\end{equation}

We know that the last term is less than the sum of all previous terms, where $t$ is sufficiently large. If every integer less than the sum of the previous terms can be found as the sum of some combination of the previous terms, one can find a partition for any integer $I$ by summing the $(t+1)^{TH}$ term $C_2^{t+1}$ and the partition for $I-C_2^{t+1}$, as $I-C_2^{t+1}<\sum_{N=1}^{t}\frac{N^2 - N + 2}{2}$. By induction, one then sees that any integer less than the sum of $\{C_2^N\}$ can be found as the sum of some combination of the terms from this sequence, if this can be done for the first $t$ terms, where $t \in \mathbb{N}$, $t > 3$. (If you can do this for the first $t$ terms, you can do it for the first $t+1$ terms, and thus the first $t+2$ terms, et cetera.) This method closely follows that used by Brown to show whether a sequence is complete \cite{one}. Here we take the first four terms of the sequence, then examine the integers that can be found as sums of the terms. 
\begin{equation}
    \{C_2^N\}=\{1, 2, 4, 7\}
\end{equation}
\begin{equation}
    1+2+4+7=14
\end{equation}
We see we have $1=1$, $2=2$, $3=1+2$, $4=4$, $5=1+4$, $6=2+4$, $7=7$, $8=1+7$, $9=7+2$, $10=1+2+7$, $11=4+7$, $12=1+4+7$, and $13=2+4+7$. This then confirms that the lazy caterer's sequence is complete, and hence also semicomplete. (Note that one can individually check cases to see that it is semicomplete when using the first 3, 2, or 1 terms. We just needed to apply this technique to a set of more than 3 terms.) To illustrate that this holds upon adding another term, we examine the first five terms of the sequence.

\begin{equation}
    \{C_2^N\}=\{1, 2, 4, 7, 11\}
\end{equation}
\begin{equation}
    1+2+4+7+11=25
\end{equation}
We already know we have partitions for $1, 2, \ldots, 14$. Now we need to find these for $15,16 \ldots, 25$. As an example, we take the case of $I=20$. We see that $20-11=9=2+7$, hence $20=11+2+7$.  

\noindent \textbf{In higher dimensions:}
There is a known recurrence relation that describes such sequences. 
\begin{equation}
    C_{M}^{N}=C_{M}^{N-1}+C_{M-1}^{N-1} 
\end{equation}
\begin{equation*}
    C_M^1=1
\end{equation*}
For a discussion of this recurrence, see \cite{five}.

\noindent \textbf{In three dimensions:}
\begin{equation}
    C_{3}^{N}=\bigg\{ \frac{N^3 - 3N^2 +8N}{6} \bigg\}
\end{equation}
\begin{equation}
    \sum_{N=1}^{t} \frac{N^3 - 3N^2 +8N}{6} = \frac{t(t+1)(t^2-3t+14)}{24}
\end{equation}
We then have
\begin{equation}
    \frac{(t+1)^3 - 3(t+1)^2 +8(t+1)}{6} < \frac{t(t+1)(t^2-3t+14)}{24}
\end{equation}

\noindent where $4 < t$.

\noindent \textbf{In four dimensions:}
\begin{equation}
    C_{4}^{N}=\bigg\{ \frac{N^4-6N^3+23N^2-18N+24}{24} \bigg\}
\end{equation}
\begin{equation}
    \sum_{N=1}^{t} \frac{N^4-6N^3+23N^2-18N+24}{24} = \frac{t(t^4-5t^3+25t^2+5t+94)}{120}
\end{equation}

\noindent We then have 
\begin{equation}
    \frac{(t+1)^4-6(t+1)^3+23(t+1)^2-18(t+1)+24}{24} < \frac{t(t^4-5t^3+25t^2+5t+94)}{120}
\end{equation}
where $5 < t$. One can then apply the same process used in two dimensions to prove semicompleteness in three and four dimensions. For three dimensions one must manually check semicompleteness for $t \leq 4$, and for four dimensions, one must manually check semicompleteness for $t \leq 5$.
\end{proof}

\section{Proof for Theorem \ref{thm3}}
\label{sec4}
\begin{proof}
\begin{equation}
    F(q)=\sum_{m=0}^{j+1} \sum_{n=0}^{j+1} \frac{(-1)^{\frac{3n}{2}}q^{\frac{1}{2}(n^2-2n+2)} (-1)^{\frac{m}{2}}q^{\frac{1}{2}(m^2-2m)}}{q+1} {j+1 \brack n}_q {j+1 \brack m}_q 
\end{equation}
We can split this double sum into the product of two sums.
\begin{equation}
    F(q)=\biggl( \sum_{m=0}^{j+1} \frac{(-1)^{\frac{3m}{2}}q^{\frac{1}{2}(m^2-2m+2)}}{q+1} {j+1 \brack m}_q \biggr) \biggl(\sum_{m=0}^{j+1} (-1)^{\frac{m}{2}}q^{\frac{1}{2}(m^2-2m)} {j+1 \brack m}_q \biggr)
\end{equation}
By noting the definition for a binomial coefficient and manipulating the series, we can arrive at a more convenient form. 
\begin{equation}
    \binom{m}{2}=\frac{m!}{2!(m-2)!}=\frac{1}{2}m(m-1)
\end{equation}
\begin{equation}
    F(q)=\biggl(\frac{q}{q+1}\biggr)\biggl(\sum_{m=0}^{j+1} (-1)^{m} \biggl(-\frac{1}{q}\biggr)^{\frac{m}{2}} q^{\binom{m}{2}} {j+1 \brack m}_q \biggr)\biggl(\sum_{m=0}^{j+1} \biggl(-\frac{1}{q}\biggr)^{\frac{m}{2}} q^{\binom{m}{2}} {j+1 \brack m}_q \biggr)
\end{equation}
We have the following identity from the work of Koekoek and Swarttouw \cite{six}. It follows from the $q$-binomial theorem.
\begin{equation}
    (a;q)_n=\sum_{k=0}^{n} (-a)^{k} q^{\binom{k}{2}} {n \brack k}_q
\end{equation}
This identity then enables us to rewrite $F(q)$.
\begin{equation}
    F(q) = \frac{q}{q+1} \biggl(\sqrt{-\frac{1}{q}};q \biggr)_{j+1} \biggl(-\sqrt{-\frac{1}{q}};q \biggr)_{j+1}
\end{equation}
The identity given below is known \cite{six}. 
\begin{equation}
    (a^2;q^2 )_{n} = (a;q )_{n} (-a;q )_{n}
\end{equation}
\begin{equation}
    \therefore \hspace{0.2 cm} F(q)=\frac{q}{q+1} \biggl(\sqrt{-\frac{1}{q}};q \biggr)_{j+1} \biggl(-\sqrt{-\frac{1}{q}};q \biggr)_{j+1} =\frac{q}{q+1} \biggl(-\frac{1}{q};q^2 \biggr)_{j+1}
\end{equation}
The $q$-analog of the Pochhammer symbol can be defined as the following product for positive $k$.
\begin{equation}
    (a;q)_n=\prod_{k=0}^{n-1} (1-aq^{k})
\end{equation}
\begin{equation}
    F(q)=\frac{q (-\frac{1}{q};q^{2})_{j+1}}{q+1} = \frac{q}{q+1} (1+q^{-1})(1+q)(1+q^{3})(1+q^{5})\cdots = \prod_{n=1}^{j} (1+q^{2n-1})
\end{equation}
This product is a generating function for partitions of $kj$ into distinct odd integers drawn from $\{1,3,5,\ldots,2j-1\}$.
\begin{equation}
    \prod_{n=1}^{j} (1+q^{2n-1}) = \sum_{kj=0}^{j^2} f_{O_j}(k,j) q^{kj}
\end{equation}

\noindent The Pell constant $\mathcal{P}_{Pell}$ can be written as shown below \cite{three}.
\begin{equation}
    \mathcal{P}_{Pell}=1-\prod_{k=0}^{\infty} \biggl(1-\frac{1}{2^{2k+1}} \biggr)
\end{equation}
\begin{equation}
    \lim_{j \to \infty} F\biggl(-\frac{1}{2}\biggr)=\lim_{j \to \infty} \biggl(\prod_{n=1}^{j} \biggl(1+\biggl(-\frac{1}{2}\biggr)^{2n-1}\biggr)\biggr)=\prod_{n=0}^{\infty} \biggl(1-\frac{1}{2^{2n+1}} \biggl) = 1-\mathcal{P}_{Pell}
\end{equation}
Hence $\mathcal{P}_{Pell}=1-\Lim{j \to \infty}(F(-\frac{1}{2}))$.
\end{proof}

\end{document}